\documentstyle[12pt]{amsart}
\input diagrams
\diagramstyle{loose,width=3em,height=2em,labelstyle=\scriptstyle}
\def\bbd{\begin{diagram}}
\def\eed{\end{diagram}}

\def\myemail{{phung@@msri.org}}
\def\myperemail{{phung@@ioit.ncst.ac.vn}}

\def\myabstract{We show that, with some technical conditions, an abelian category can be embedded into the category of bimodules over a ring. The case of semisimple rigid monoidal categories is studied in more detail.}
\def\myaddress{Mathematical Sciences Research Institute, 1000 Centennial Drive, Berkeley, CA 94720}
\def\myperaddress{Hanoi Institute of Mathematics, P.O.Box 631, 10000 Boho, Hanoi}
\def\mythanks{The author should like to thank Professors M. Artin, B. Pareigis and P. Smith for useful discussions. This work is supported by the Mathematical Sciences Research Institute.

This work was presented at the conference on Quantum groups in Morelia, Mexico. The author should like to thank the organisers and the Institute of Mathematics, UNAM for the hospitality during the conference.}

\def\myauthor{Ph\`ung H{{\accent"5E o}\kern-.28em\raise.2ex\hbox{\char'22}\kern-.20em} H{a\kern-.370em\raise.16ex\hbox{\char'47}\kern.1em}i}
\def\myAUTHOR{ PH\`UNG H{{\accent"5E O}\kern-.38em\raise.8ex\hbox{\char'22}\kern-.12em}  H{A\kern-.46em\raise.80ex\hbox{\char'47}\kern.18em}I}

\def\amshead{
\title{An Embedding theorem for abelian monoidal categories}
\author{ PH\`UNG H{{\accent"5E O}\kern-.38em\raise.8ex\hbox{\char'22}\kern-.12em}  H{A\kern-.46em\raise.80ex\hbox{\char'47}\kern.18em}i}
\address{\myperaddress}
\curraddr{\myaddress}
\keywords{embedding theorem, monoidal category}
\thanks{{\it 2000 Mathematics Subject Classification.} Primary 18D10, 18E20, Secondary 16D20, 16W30}
\email{\myperemail \ {}{\rm and} \myemail}
\begin{abstract}\myabstract
\end{abstract}
\maketitle }

\def\lora{\longrightarrow}
\def\lola{\longleftarrow}
\def\ot{\otimes}

\def\loma{\longmapsto}

\def\alph{\alpha}

\newcommand{\bbas}{\begin{eqnarray*}}
\newcommand{\eeas}{\end{eqnarray*}}
\newcommand{\bbar}{\begin{array}}
\newcommand{\eear}{\end{array}}
\newcommand{\bbs}{\begin{displaymath}}
\newcommand{\ees}{\end{displaymath}}
\newcommand{\bb}{\begin{equation}}
\newcommand{\eqbb}{\begin{equation}}
\def\ee{\end{equation}}
\def\eqee{\end{equation}}
\def\eea{\end{eqnarray}}

\def\bba{\begin{eqnarray}}
\newtheorem{thm}{Theorem}[section]
\newtheorem{lem}[thm]{Lemma}

\newtheorem{cor}[thm]{Corollary}

\newtheorem{pro}[thm]{Proposition}

\def\Im{\mbox{\sf Im}}

\def\Hom{\mbox{\sf Hom}}

\def\End{\mbox{\sf End}}

\def\V{{\cal V}}

\def\eee{\rule{.75ex}{1.5ex}\\[1ex]}

\def\proof{{\it Proof.\ }}

\renewcommand{\dim}{\mbox{\sf dim}}

\def\rref#1{(\ref{#1})}

\font\Fraktur=eufm10 scaled\magstep1          
   \newcommand{\fraktur}[1]{\mbox{\Fraktur #1}}  %
   \font\Fraktu=eufm7 scaled\magstep1            
   \newcommand{\fraktu}[1]{\mbox{\Fraktu #1}}    %
   \font\Frakt=eufm5 scaled\magstep1             
  \newcommand{\frakt}[1]{\mbox{\Frakt #1}}      %
   \def\frak#1{\mathchoice{\fraktur {#1}}            
                        {\fraktur {#1}}            
                        {\fraktu {#1}}             
                        {\frakt {#1}}  }           

\newcommand{\Ss}{\frak S}

\newcommand{\bK}{{K}}
\newcommand{\bZ}{{\Bbb Z}}

\def\db{{\mathchoice{\mbox{\sf db}}
                    {\mbox{\sf db}}
                    {\mbox{\scriptsize\sf db}}
                    {\mbox{\tiny\sf db}} }}

\def\ev{{\mathchoice{\mbox{\sf ev}}
                    {\mbox{\sf ev}}
                    {\mbox{\scriptsize\sf ev}}
                    {\mbox{\tiny\sf ev}} }}

\def\id{{\mathchoice{\mbox{\sf id}}
                    {\mbox{\sf id}}
                    {\mbox{\scriptsize\sf id}}
                    {\mbox{\tiny\sf id}} }}

\def\op{{\mathchoice{\mbox{\rm op}}
                    {\mbox{\rm op}}
                    {\mbox{\scriptsize\rm op}}
                    {\mbox{\tiny\rm op}} }}

\newcommand{\Mod}{{\sf Mod}}

\def\part{\vdash}

\def\lam{{\lambda}}

\def\Fun{{\sf Fun}}
\def\ind{{\sf ind}}
\def\rhom{{\sf rhom}}
\def\lhom{{\sf lhom}}
\def\cat#1{{\text{\sf #1}}}
\def\fun#1{{\text{\sf #1}}}
\def\bigcirc{\Box}

\let\stackrell=\stackrel
\def\stackrel#1#2{\stackrell{\scriptstyle #1}{\scriptscriptstyle #2}}

\let \liml=\lim
\def\lim{\displaystyle\liml}
\def\V#1#2#3{{V_{#1}\ }^{#2}_{#3}}

\begin{document}
\bibliographystyle{plain}
\amshead
\section*{Introduction}
The problem of finding an embedding theorem for abelian monoidal categories is motivated on one hand by the Freyd-Mitchell full embedding theorem and on the other hand by Deligne's and Doplicher-Roberts' theories of ``abstract Tannakian-Krein's duality'' \cite{deligne90,dr89}.

By definition, an embedding from an abelian category to another is a faithful functor which sends non-zero objects to non-zero objects.
According to Freyd and Mitchell, a small abelian category admits an exact full embedding into the category of modules over a ring \cite{mitchell}. This important theorem allows one for example to treat finite diagrams in an abelian category as diagrams of modules.

In \cite{deligne90}, Deligne shows that, under certain technical conditions, an abelian symmetric rigid monoidal category admits an exact monoidal embedding into the category of modules over a commutative ring. Hence, by Tannaka-Krein's duality, such a category is equivalent to the category of representation of a groupoid \cite[Theorem~1.12]{deligne90}. Doplicher and Roberts study the case of compact groups and obtain an analogous result for $C^*$-categories \cite{dr89}.

With the birth of quantum groups \cite{drinfeld86}, the theory of monoidal category has a new motivation. The ``symmetric'' condition turns out to be too strong and is replaced by a weaker one, the ``braided'' condition. The problem of generalizing Deligne's and Doplicher-Roberts' results to braided categories is interesting. Yet, one does not know what can be a target for such an embedding.

The situation seems to be simpler if we drop the symmetry, i.e., to find an embedding for abelian monoidal categories. A natural candidate for the target category is the category of bimodules over a ring. In this paper we show that any small monoidal cateory with exact tensor product admits a right exact monoidal embedding into the category of bimodules over a ring. In particular, a small abelian rigid monoidal category admits an exact monoidal embedding (Theorem \ref{thm2}).

Unfortunately, this embedding theorem does not seems to help solve the problem for braided categories, for, according to Schauenburg \cite{schauen1}, the category of bimodules over a ring is never braided unless the basic ring is a field. 

Let me briefly explain the main idea of the paper. To find an embedding for a small abelian monidal category, we first extend it to a bigger monoidal category which is cocomplete and has an injective cogenerator, namely, a Grothendieck monoidal category. Then we extend the latter category to a module category, say over a ring $R$. Finally, we construct a monoidal functor from a module category with a monoidal structure to a bimodule category with the usual tensor product. The construction of the last functor will be given in Section \ref{sec2}. In Section \ref{sec3} we explain how to extend the tensor product on a small abelian monoidal category to a monoidal structure on a module category containing this small abelian monoidal category. In the last section we consider an application to the special case of semisimple categories. An explicit embedding is described. As a consequence, we show that a semisimple symmetric monoidal category with simple unit object is Tannakian (Corollary \ref{cor2}).

Throughout the paper, the tensor product over a ring $R$ is denoted by $\ot$. $\ot$ also denotes the tensor product in an abstract monoidal category when no confusion may appear, othewise, we use the signs $\odot$ or $\Box$. The category of right $R$-modules (resp. left $R$-modules or $R-R$-bimodules) is denoted by $\Mod_R$ (resp. ${}_R\Mod$ or ${}_R\Mod_R$). $\Hom_R$ (resp. ${}_R\Hom$ or ${}_R\Hom_R$) denotes the set of homomorphisms of right $R$-modules (resp. left $R$-modules or $R-R$-bimodules).

\section{Abelian Monoidal categories}\label{sec1}
\subsection{Monoidal categories}\label{sec1.1}

Let $\cat A$ be a category. A monoidal structure on $\cat A$ consists of the following data: a bifunctor $\ot:\cat A\times \cat A\lora \cat A$, $(X,Y)\loma X\ot Y$, called tensor product, an object $I$, called unit object, for which
\begin{enumerate}\item there exists a natural isomorphism $\alpha$ between the functors $(-\ot -)\ot -$ and $-\ot(-\ot -)$, $\alpha_{X,Y,Z}:(X\ot Y)\ot Z\lora X\ot (Y\ot Z)$, called associator, such that the following diagram commutes:
\bba\label{eq1}
\begin{diagram}
((X\ot Y)\ot Z)\ot U &\rTo^{\alpha_{X,Y,Z}\ot 1} &(X\ot (Y\ot Z))\ot U &\rTo{\alpha_{A,B\ot C,D}}& X\ot(( Y\ot Z)\ot U)\\
\dTo<{\alpha_{X\ot Y,Z,U}}&&&& \dTo>{1\ot\alpha_{Y,Z,U}}\\
(X\ot Y)\ot (Z\ot U)&&\rTo^{\alpha_{X,Y,Z\ot U}}&& X\ot (Y\ot (Z\ot U))
\end{diagram}
\eea
\item there exist natural ismorphisms of functors $-\ot I$, $I\ot -$ and the identity functor: $\rho_X:X\ot I\lora X$ and $\lam_X:I\ot X\lora X$, called right and left units, such that the following diagram commutes:
\bba\label{eq2}\begin{diagram}
(X\ot I)\ot Y&&\rTo^{\alpha_{X,I,Y}}&&X\ot (I\ot Y)\\
&\rdTo<{\rho_X\ot 1}&&\ldTo>{1\ot\lam_Y}\\
&&X\ot Y&&\end{diagram}\eea
\end{enumerate}

$(\cat A,\ot, I,\alpha,\lam,\rho)$ is called a monoidal category. In the case, when $\alpha,\lam,\rho$ are identity morphisms, we have a strict monoidal category. In the general case, the associator $\alph$ allows one to speak of a tensor product of many objects $X_1\ot X_2\ot\cdots\ot X_n$, without specifying the order in which the tensor product is applied (MacLane's coherence theorem \cite[VII,2]{maclane}).

Using \rref{eq1} and \ref{eq2}, we can show that the following diagrams commute:
\bba\label{eq2.1}
\begin{diagram}
(I\ot X)\ot Y&&\rTo^{\alpha_{I,X,Y}}&&X\ot (I\ot Y)\\
&\rdTo<{\lam_X\ot 1}&&\ldTo>{\lam_{X\ot Y}}\\
&&X\ot Y&&\end{diagram}\qquad\begin{diagram}
(X\ot Y)\ot I&&\rTo^{\alpha_{X,Y,I}}&&X\ot (Y\ot I)\\
&\rdTo<{\rho_{X\ot Y}}&&\ldTo>{1\ot\rho_Y}\\
&&X\ot Y&&\end{diagram}\eea
and that $\lam_I=\rho_I$ (cf. \cite{par3}).

From the definition of a bifunctor, we have
\bba\label{eq3} f\ot g=(1\ot g)\circ(f\ot 1)=(f\ot 1)\circ(1\ot g).\eea
In particular, using the isomorphism $\lam_I=\rho_I$, we have, for $r,s\in \End(I)$:
\bba\label{eq4}\begin{diagram}
I\ot I&\rTo^{r\ot 1}&I\ot I&\rTo^{s\ot 1}& I\ot I\\
\dTo>{\rho_I}&&\dTo<{\rho_I}>{\lam_I}&&\dTo<{\lam_I}\\ I&\rTo^{r}&I&\rTo^{s}&I.\end{diagram}\eea
Therefore, using the equation \rref{eq3}, we have $r\circ s=s\circ r$. Thus, $\End(I)$ is an abelian group with respect to composition. Further, this group acts on any set $\Hom(X,Y)$ from the left and the right by means of the isomorphism $\lam$ and $\rho$:
\bba\label{eq5} r\cdot f:=\lam_Y(r\ot f)\lam_X^{-1},\quad f\cdot r:=\rho_Y(f\ot r)\rho_X^{-1}.\eea
We have $1_I\cdot f=f\cdot 1_I=f$.

\subsection{The internal homs}\label{sec1.2}
Each object $X$ in $\cat A$ defines a functor $X\ot -:\cat A\lora \cat A$, $Y\loma X\ot Y$. If this functor has a right adjoint, the right adjoint will be denoted by $\rhom(X,-)$. We have, by definition, a natural isomorphism
\bba\label{eq7} \Hom(X\ot Y,Z)\cong\Hom(Y,\rhom(X,Z)), \quad\forall\ Y,Z.\eea
The functor $\lhom$ is define analogously by 
\bba\label{eq8} \Hom(Y\ot X, Z)\cong \Hom(Y,\lhom(X,Z)),\quad \forall\ Y,Z.\eea

The category $\cat A$ is called left closed (resp. right closed or closed) if the functor $\lhom(X,-)$ (resp. $\rhom(X,-)$ or both functors)  is defined for any $X\in \cat A$. The general theory of adjoint functors (cf. \cite[Corollary V.3.2]{mitchell}) gives us the following criteria of closedness
\begin{lem}\label{lem1} Assume that $\cat A$ is a cocomplete category with a generator. Then the tensor product on $\cat A$ is closed if and only if it commutes with colimits.\end{lem}

In the general case, the functor $\lhom(X,-)$ and $\rhom(X,-)$ preserves colimits, whenever the latter is defined. If the category is closed, then we can also speak of the (contravariant) functors $\lhom(-,Z)$ and $\rhom(-,Z)$.
\begin{lem}\label{lem2}
Assume that $\cat A$ is closed. Then the functors $\lhom(-,Z)$ and $\rhom(-,Z)$ preserve colimits (i.e. sending colimits to limits) whenever the latter is defined.\end{lem}
\proof We have, for any $Y\in\cat A$,
\bbas \Hom(Y,\rhom(\lim_{\lora} X_i,Z)&\cong& 
\Hom(\lim_{\lora} X_i\ot Y,Z)\\
&\cong& \Hom(\lim_{\lora}(X_i\ot Y),Z)\\
&\cong&\lim_{\longleftarrow}\Hom(X_i\ot Y,Z)\\
&\cong&\lim_{\longleftarrow}\Hom(Y,\rhom(X_i,Z))\\
&\cong&\Hom(Y,\lim_{\longleftarrow}\rhom(X_i,Z)),\eeas
here we use the fact that $\Hom(Y,-)$ preserves limits. Since the above isomorphisms hold for any $Y$, we conclude that
\bbs \rhom(\lim_{\lora} X_i,Z)\cong\lim_{\longleftarrow}\rhom(X_i,Z).\ees

The case of $\lhom(-,Z)$ is treated quite analogously.\eee

\subsection{Rigid objects}\label{sec1.3}

Setting $Z=I$ and $Y=\lhom(X,I)$ in \rref{eq8}, we obtain a morphism $\ev_X:\lhom(X,I)\ot X\lora I$, which corresponds to the identity morphism in $\Hom(\lhom(X,I),\lhom(X,I))$. It is obvious that $(\lhom(X,I),\ev_X)$ is universal with this property. By definition, a left dual to an object $X$ is a pair $(X^*,\ev_X:X^*\ot X\lora I)$ such that there exists a morphism $\db_X:I\lora X\ot X^*$ making the following diagrams commutative:
\bba\label{eq9}\begin{diagram}
X&\rTo<{\db_X\ot 1}&(X\ot X^*)\ot X\\
\dTo>{\id}&&\dTo>{\alpha_{X,X^*,X}}\\
X&\lTo^{\id_X\ot\ev_X}&X\ot(X^*\ot X)\end{diagram}\qquad\qquad
\begin{diagram}
X^*&\rTo<{\id\ot \db_X}&X^*\ot(X\ot X^*)\\
\dTo>{\id_{X^*}}&&\dTo>{\alpha_{X^*,X,X^*}^{-1}}\\
X&\lTo^{\ev_{X}\ot\id_{X^*}}&(X^*\ot X)\ot X^*\end{diagram}
\eea
The left dual, if it exists, is uniquely determined up to an isomorphism. Moreover, we have a natural isomorphism
\bba\label{eq10}
\lhom(X,Z)\cong Z\ot X^*.\eea
In particular, the functor $\lhom(X,-)$ exists if $X$ has a left dual. In this case, the diagrams in \rref{eq9} also imply that the functor $X\ot -$ has a left adjoint: $X^*\ot -$, consequently $X\ot -$ commutes with limits.

The definition of a right dual ${}^*X$ to $X$ is similar.
 Analogous assertions hold for objects having right dual. An object in $\cat A$ is called rigid if it possesses left and right duals. The category $\cat A$ is called rigid if its objects are rigid.

\subsection{Monoidal functors}\label{sec1.4}
Let $(\cat A,\ot,\alpha,\rho,\lam)$ and $(\cat A',\ot',\alpha',\rho',\lam')$ be monoidal categories. A functor $\fun F:\cat A\lora \cat A'$ is called monoidal functor if there exists a natural ismorphism 
$\xi_{X,Y}:\fun F(X)\ot'\fun(Y)\lora \fun F(X\ot Y)$ and an isomorphism $\eta:I'\lora \fun F(I)$,
such that the following diagrams are commutative:
\bba\label{eq6}\begin{array}{c}\begin{diagram}
(\fun F(X)\ot' \fun F(Y))\ot' \fun F(Z)&\rTo^{\xi\ot' 1}&\fun F(X\ot Y)\ot'\fun F(Z)&\rTo^{\xi}&
\fun F((X\ot Y)\ot Z)\\
\dTo>{\alpha'}&&&&\dTo<{\fun F(\alpha)}\\
\fun F(X)\ot'(\fun F(Y)\ot'\fun F(Z))&\rTo^{1\ot' \xi}&
\fun F(X)\ot' \fun F(Y\ot Z)&\rTo^{\xi}&\fun F(X\ot (Y\ot Z))\end{diagram}\\[4ex] \\ 
\begin{diagram}
I'\ot'\fun F(X)&\rTo^{\eta\ot 1}&\fun F(I)\ot'\fun F(X)\\
\dTo>{\lam'}&&\dTo<{\xi}\\
\fun F(X)&\lTo^{\lam'}&\fun F(I\ot X)\end{diagram}\qquad
\begin{diagram}
\fun F(X)\ot' I'&\rTo^{1\ot \eta}&\fun F(X)\ot'\fun F(I)\\
\dTo>{\rho'}&&\dTo<{\xi}\\
\fun F(X)&\lTo^{\fun F(\rho)}&\fun F(X\ot I)\end{diagram}\eear\eea

If an object $X$ from  $\cat A$ is rigid then it image $\fun F(X)$ is also rigid.
\subsection{Abelian monoidal categories}\label{sec1.5}
A monoidal category $(\cat A,\ot)$ is called abelian monoidal if it is abelian and the tensor product is an additive bifunctor.

In this case, $\bK:=\End(I)$ is a commutative ring and $\Hom(X,Y)$ becomes $\bK-\bK$-bimodule, for any objects $X,Y$. Notice that the two actions of $\bK$ do not generally coincide and $\cat A$ is therefore not necessarily $\bK$-linear. 

From the discussion in \ref{sec1.2}, if the functor $\lhom(X,-)$ is defined then it is left exact and the functor $X\ot -$ is right exact. The same holds for $\rhom(X,-)$ and $-\ot X$. If $\cat A$ is closed then the contravariant functor $\lhom(-,Z)$ and $\rhom(-,Z)$ are also left exact. If an object $X$ has left dual then the functor $X\ot -$ is left exact; hence exact. In particular, if $X$ is rigid then all the mentioned above functors are exact. Since $X^*=\lhom(X,I)$, for a short exact sequence of rigid objects $0\lora X\lora Y\lora Z\lora 0$ its dual sequence $0\lora Z^*\lora Y^*\lora X^*\lora 0$ is also exact.

\section{Tensor structures on module categories}\label{sec2}

Let $R$ be a a ring and $\Mod_R$ be the category of right modules over $R$. Let $\odot$ be a tensor product on $\Mod_R$ with associator $\alpha$ and left-, right units $\lam,\rho$. We assume that the tensor product is (left and right) closed. That is, (cf \ref{lem1}) the functors $-\odot X$ and $X\odot -$ possess right adjoints. Consequently, these functors preserves colimits.

Consider a functor $-\odot X:\Mod_R\lora \Mod_R$, which preverves colimits. By a theorem of Watts \cite{watts}, there exists a left action of $R$ on $R\odot X$ making it an $R-R$-bimodule, such that the functor $-\odot X$ is naturally equivalent to the functor $-\ot (R\odot X)$: $Y\loma Y\ot (R\odot X)$. That is, there is a natural isomorphism
\bba\label{eq2.2}\theta_X(Y):Y\ot (R\odot X)\lora Y\odot X.\eea
Throughout this section, $\ot$ denotes the tensor product over $R$.
For the case $X=R$, we shall call the action of $R$ on $T:=R\odot R$ the first left action of $R$ on $T$, to distinguish with the second action defined subsequently.

Explicitly, the left action of $R$ on $R\odot X$ is given as follows
\bba\label{eq2.3}
-\odot X:R\cong\End_R(R)\lora \End_R((R\odot X)_.).\eea
Since $\odot$ is biadditive, we see that $\theta_X(Y)$ is also natural on $X$. Hence
$R\odot -:X\loma R\odot X$ is a functor form the category $\Mod_R$ to $_R\Mod_R$, commuting with colimits. Applying Watts' theorem again, we have an equivalence of the functors $-\odot R$ and $-\ot  T$
with a left action of $R$ on $T$, called the second left action:
\bba\label{eq2.4} \mu_X:R\odot X\lora X\ot_2T,\eea
where the subindex 2 indicates that the second left action of $R$ on $T$ is used to define the tensor product.

By its definition, the second left action commutes with the other actions of $R$ on $T$, making $T$ an object in $_{(R\ot_\bZ R)}\Mod_R$. And have an $R$-linear natural isomorphism
\bba\label{eq2.5} c_{X,Y}:=\theta_Y(X)\circ(\id_X\ot\mu_Y):X\ot  {}_2(Y\ot  {}_1T)\lora X\odot Y.\eea

The associator $\alph$ induces an $R$-linear natural isomorphism
\bba\nonumber
\alpha'_{X,Y,Z}&:=&(\id\ot c_{Y,Z})c^{-1}_{X,Y\odot Z}\alpha_{X,Y,Z}c_{X\odot Y,Z}(c_{X,Y}\ot \id):\\
\label{eq2.5a}
&& (X\ot_2(Y\ot_1T))\ot_2(Z\ot_1T_\cdot)\lora X\ot_2((Y\ot_2(Z\ot_1T))\ot_1T_\cdot)\eea
where the ation of $R$ is induced from the right action on $T$ (indicated by a dot).
Analogously, we have $R$-linear natural isomorphisms $\lam'$ and $\rho'$
\bba\label{eq2.6} \lam'_X&:=&\lam_X\circ \xi_{I,X}:I\ot_2(X\ot_1T)\lora X,\\
\label{eq2.6a} \rho'_X&:=&\rho_X\circ \xi_{X,I}:X\ot_2(I\ot_1T)\lora X.\eea
Thus, we have defined data for a monoidal structure on $\Mod_R$, namely, the tensor product of two module $X,Y$ is $X\ot_2(Y\ot_1T))$, the associator is $\alpha'$ and the left and right units is $\lam'$ and $\rho'$. Using the naturality of $c$, we can show that these data define a monoidal structure on $\Mod_R$. More explicitly, using routine diagram chasing, we can show the following lemma.
\begin{lem}\label{lem3} Let $(\cat A,\ot,I,\alpha,\rho,\lam)$ be a monoidal category. Let $\odot$ be another bifunctor $\cat A\times \cat A\lora \cat A$, which is equivalent to $\ot$ by means of a natural isomorphism $\xi$. Then $\odot$ together with the isomorphism $\alpha'$, $\rho'$ and $\lam'$ defined as in \rref{eq2.5a},\rref{eq2.6} and \rref{eq2.6a}, define another monoidal structure on $\cat A$.\eee\end{lem}

Setting $X=Y=Z=R$ in \rref{eq2.5a}, we obtain an isomorphism of $R\ot_\bZ R\ot_\bZ R-R$-bimodules
\bba\label{eq2.7}\alpha_{R,R,R}:{}_{{\circ\atop {\uparrow\atop 1}}{\bullet\atop {\uparrow\atop 2}}}T\ot_2{}_{*\atop {\uparrow\atop 3}}T_\cdot\lora {}_{{*\atop {\uparrow\atop 3}}{\circ\atop {\uparrow\atop 1}}}T\ot_1{}_{\bullet\atop {\uparrow\atop 2}} T_\cdot\eea
where $\bullet,\circ,*$ denote the different actions of $R$ on the source and the target of $\Psi$ they will be referred to as the first, second and third left action of $R$. The right action is indicated by $\cdot$. The following lemma shows that $\alpha'$ can be restored from this isomorphism.
\begin{lem}\label{lem4} Let $P,Q$ be $R-S$-bimodules. Then a natural transformation $c:-\ot P\lora -\ot Q$ of functors $\Mod_R\lora\Mod_S$ is given by a $R-S$-bimodule homomorphism $c=c_R:P\lora Q$.\end{lem}
\proof We have the following commutative diagram:
\bbs\begin{diagram}[loose,width=3em,height=2em]
M\ot P&\rTo^{c_M}&M\ot Q\\
\dTo<{f}&&\dTo>{f\ot\id_Q}\\
N\ot P&\rTo^{c_N}&N\ot Q\end{diagram}\ees
For $N=M=R$ and a morphism $f_s:R\lora R$ or right $R$-modules, $f_s(r):=sr$, we have
\bbs c(sp)=c_R(s\ot p)=s(c_R(1\ot p))=sc(p).\ees That is $c$ is a left $R$-module morphism. By definition, $c$ is a right $S$-module morphism, hence it is an $R-S$-bimodule morphism. Let now $M=R$ and $N$ arbitrary. For $n\in N$, choose $f_n:R\lora N$, $f_n(s)=ns$, thus $f$ is a morphism of right $R$-modules. Then we have, plugging $f$ and $c_N$ in the above diagram, $c_N(n\ot p)=(f\ot\id_Q)c_R(1\ot p)=n\ot c(p).$\eee

Setting in \rref{eq2.6a} $X=R$ and $Y=I$, we have a natural isomorphism
\bba\label{eq2.8}
\rho':=\rho\circ c(R,I):I\ot_1T\lora R.\eea
Analogously, we have a natural isomorphism
\bba\label{eq2.9}
\lam':=\lam\circ c(I,R):I\ot_2T\lora R.\eea

The coherent constraints \rref{eq1},\rref{eq2}, imply the following condition for the isomorphisms $\alpha',\lam',\rho'$:
\bba\label{eq2.9a} (\alpha'\ot_1\id_T)(\id_T\ot_1\alpha')(\alpha'\ot_2\id_T)&=& (\id_T\ot_3\alpha')(\id_T\ot_2\alpha'),\\
\label{eq2.9b}
(\lam'\ot_1\id_T)(\rho'\ot_2\id_T)&=&\id_I\ot_1\alpha'.\eea

For any right $R$-modules $M,N$, we have a sequence
\bba\label{eq2.10}
\bbd \Hom_R(M,N)&\rTo^{\id_M\ot -} &_2\Hom_{R-R}(M\ot_1T_\cdot,N\ot_1T_\cdot)&\rTo^{\id_I\ot -} &  \Hom_R(M_\cdot,N_\cdot)\\
f&\rMapsto&f\ot \id_T&\rMapsto&\id_I\ot f\ot \id_T= f\eed\eea
here $_2\Hom_R(M\ot_1T_\cdot,N\ot_1T_\cdot)$ denotes the set of $R-R$-bimodule morphisms, the left actions of $R$ on whose source and target are given by the second left action of $R$ on $T$. According to the isomorphism in \rref{eq2.6} the composition of the above morphisms is an isomorphism, hence the map
\bba \bbd\Hom_R(M,N)&\rTo& {}_2\Hom_R(M\ot_1T,N\ot_1T)\eed\eea
is injective.

Consider now the functor $\omega:\Mod_R\lora \Mod_R$, $X\loma X\ot_1T$, where the left action of $R$ on $X\ot_1T$ is given by the second left action of $R$ on $T$. We have a natural isomorphism
\bba\nonumber
{\xi_{X,Y}:=\alpha'_{R,X,Y}}:&&\omega(X)\ot\omega(Y)=(X\ot_1T)\ot_2(Y\ot_1T)
\\
\label{eq2.11}
&&\stackrel{\alpha_{R,X,Y}}{\textstyle\lora}(X\ot_2T))\ot_1T=\omega(X\ot_2(Y\ot_1T)).\eea
\begin{lem}\label{lem5} With the notation as above, $(\omega,\xi)$ is a monoidal functor from $(\Mod_R,\odot)$ to $(_R\Mod_R,\ot)$.\end{lem}
\proof One has to check hexagon identity for $\xi$. Since in $_R\Mod_R$, the tensor product is strict, the hexagon can be reduced to a pentagon. It turns out that the commutativity of this pentagon is precisely the coherence condition of $\alpha'$, which holds by Lemma \rref{lem3}.\eee

\noindent{\it Remark.} Although our monoidal category $(\Mod_R,\odot)$ is not assumed to be strict. The functor $\omega$ maps it into a strict monoidal category.

\begin{thm}\label{thm1} Let $\odot$ be a monoidal structure on $\Mod_R$. Then the functor $\omega:\Mod_R\lora {_R\Mod_R}$, $X\loma X\ot_1T$ is a monoidal right exact embedding.\end{thm}
\proof We have seen that $\omega$ is monoidal, faithful, right exact. It remains to show that $\omega(X)\not \cong 0$  whenever $X\not\cong 0$. That is, $X\odot R\not\cong 0$ for all $X\not\cong 0$. Let $X$ be such that $X\odot R\cong 0$. Then 
\bbs 0=\Hom(X\odot R,X)\cong\Hom(R,\rhom(X,X)).\ees
Therefore $\rhom(X,X)= 0$. But then $\Hom(X,X)\cong\Hom(I,\rhom(X,X))=0$. Thus, $X\cong 0$.\eee

Let now $\cat A$ be a monoidal category which, as an abelian category, is cocomplete with a progenerator (i.e. small projective generator \cite{mitchell}), and, as a monoidal category, is closed. Let $P$ be a progenerator of $\cat A$. Then $\cat A$ is equivalent to $\Mod_R$, where $R=\End(P)$, by the functor $\fun F=\Hom(P,-)$. Since $\fun F$ is an equivalence, it caries the monoidal structure on $\cat A$ over to $\Mod_R$. Thus, we have
\begin{cor}\label{cor1}
Let $(\cat A,\odot)$ be an abelian monoidal category, cocomplete and closed and with a progenerator. Then the functor $\omega:X\loma \Hom(P,X\odot P)$ is a monoidal  functor from $\cat A$ to $_R\Mod_R$, where $R=\End(P)$. This functor is a right exact embedding. It is exact if $P$ is left flat with respect to the tensor product on $\cat A$. An analogous assertion holds for the functor $X\loma \Hom(P,P\odot X)$.\end{cor}

\noindent{\it Remark.} Since $X\ot_1T\cong X\odot R$. The functor $\omega$ is exact if and only if $R$ is flat with respect to the tensor product $\odot$. This is not alway the case. Take for example the category $_S\Mod_S$, where $S$ is a ring not flat over $\bZ$. Then for $R=S^{\rm op}\ot_\bZ S$, $_S\Mod_S$ is equivalent to $\Mod_R$. The tensor product is taken over $S$, therefore $M\ot_S R=M\ot_S(S^{\rm op}\ot_\bZ S)\cong M\ot_\bZ S$. Thus $R$ is not flat.

\section{An Embedding theorem for small abelian monodial categories}\label{sec3}
Using the result of the previous section, we show in this section that a small abelian monoidal category with exact tensor product can be embedded in the category of bimodules over a ring. Our tactic is to embed $\cat C$ in a bigger category which is cocomplete with a projective generator.

\subsection{The category $\cat{Ind-C}$}\label{sec2.1}
A category $\cat I$ is called a filtering category if to every pair $i,i'$ of objects from $\cat I$ there exists an object $i''$, such that $\Hom(i,i'')$ and $\Hom(i',i'')$ are both not empty, and for every pair of morphisms $f,f':i\lora i'$, there exists a morphism $g:i'\lora j$ equalizing them, i.e. $gf=gf'$.

Let $C$ be an abelian category. The category $\cat{Ind-C}$ of ind-objects of $\cat C$ consists of functors $\fun X:\cat I\lora \cat C$, where $\cat I$ is any small filtering category. Alternatively, denoting $X_i:=\fun X(i), i\in\cat I$, an ind-object of $\cat C$ is a directed system indexed by a small filtering category $\cat I$. For two objects $\fun X=\{X_i\}_{i\in\cat I}$ and $\cat Y=\{Y_j\}_{j\in\cat J}$, their hom-set is
\bba\label{eq3.1} \Hom(\fun X,\fun Y):=\lim_{\stackrel\lola i}(\lim_{\stackrel\lora  j}\Hom(X_i,Y_j)).\eea
The following lemma will be usefull when dealing with hom-sets of ind-objects.
\begin{lem}\label{lem6}(cf \cite[Appendix,~Cor. 3.2]{artin}) A morphism $f:\fun X\lora \fun Y$ can be represented, up to isomorphism, by a small filtering system of morphisms $\{f_i:X_i\lora Y_i\}_{i\in\cat I}$.\end{lem}

$\cat C$ is fully  embedded in $\cat{Ind-C}$ by a constant functor. On the other hand, $\cat{Ind-C}$ is fully embedded in $\Fun(\cat C^{\rm op},\cat{Set})$. For an ind-object $\cat X=\{X_i\}_{i\in\cat I}$, define the functor
\bba\label{eq3.2}\fun  L_{\cat X}:Y\lora \lim_{\lora}\Hom(Y,X_i).\eea
By Yoneda's Lemma and the fact that $\Hom(-\ot X)$ commutes with colimits, we have $\Hom(\cat X,\cat Y)\cong\Hom(\fun L_{\cat X},\fun L_{\cat Y})$. In fact, from definition, $\fun L_X$ is isomorphic to $\lim_{\stackrel\lora  i}\fun L_{X_i}$ in $\Fun(\cat C^{\rm op},\cat{Set})$. Therefore
\bbas \Hom(\fun L_X,\fun L_Y)&=& \lim_{\stackrel\lola i}\Hom(\fun L_{X_i},\fun L_Y)\\
(\mbox{by Yoneda's Lemma})&=&\lim_{\stackrel\lola i}\fun L_Y(X_i)\\
&=&\lim_{\stackrel\lola i}(\lim_{\stackrel\lora  j}\Hom(X_i,Y_j)).\eeas

The category $\cat{Ind-C}$ is closed under filtering direct limits (cf. \cite[Appendix. 4.4]{artin} or \cite[I.8]{sga4}).

Note however that the direct limit computed in $\cat C$ (if it exists) is generally different from the one computed in $\cat{Ind-C}$.

\subsection{Extension of functors.} \label{sec2.2}
Let $F:\cat C\lora \cat D$ be a functor. A functor $\fun{ind-F}:\cat{Ind-C}\lora \cat{ind-D}$ is defined as follows:
\bba\label{eq3.3} \fun{ind-F}(\{X_i\}_{i\in\cat I}):=\{F(X_i)\}_{i\in I}.\eea
The action of $\fun{ind-F}$ on hom-sets is defined in a straightforward manner.

\subsection{ Ind-category for abelian categories.} \label{sec2.3}
Assume now that $\cat C$ is abelian. Then the functor $\fun L_{\fun X}$ is exact for any ind-object $\fun X=\{X_i\}_{i\in\cat I}$. Indeed, in the category of sets, the filtering direct limits preserves left exact sequences, hence  for  a left exact sequence in $\cat C^{\rm op}:$ $0\lora Y\lora Y'\lora Y''$, we have the following left exact sequences
\bbs 0\lora \Hom(Y,X_i)\lora \Hom(Y',X_i)\lora \Hom(Y'',X_i)\ees
$i\in \cat I$. Taking limit we have
\bbs 0\lora \fun L_{\fun X}(Y)\lora \fun L_{\fun X}(Y')\lora \fun L_{\fun X}(Y').\ees

Conversely, let $\fun L$ be a functor $\cat C^{\rm op}\lora\cat{Ab}$. By Yoneda's Lemma, for any $X\in \cat C$, 
$\Hom(\Hom(-\ot X),\fun L(-))\cong \fun L(X)$. Consider the system $C_{\fun L}:=\{(X,\eta)|X\in \cat C,\eta\in\fun L(X)\}$. Morphism $(X,\eta)\lora (X',\eta')$ are those from $\Hom(X,X')$, which commute with $\eta,\eta'$, in the sense that $\fun L(f)(\eta')=\eta$. Then $\fun L$ is isomorphic to the functor  
$$Y\loma \lim_{\stackrel\lora  {(X,\eta)\in C_L}}\Hom(Y,X).$$ 
If $\fun L$ is left exact then $C_{\fun L}$ is a filtering system. In fact, given $(X,\eta)$ and $(X',\eta')$, since left exact functors preserve direct sums, we can take the object $(X'',\eta'')$ to be $(X\oplus X',\eta\oplus\eta')$. For any two morphisms $f,f':(X,\eta)\lora(X',\eta')$, let $g:X'\lora Y$ be the coequalizer of $f$ and $f'$. Since $\fun F$ is left exact (from $\cat C^{\rm op}$ to $\cat{Ab}$), we can take $\zeta:=\fun L(g)^{-1}(\eta')$. Then $(Y,\zeta)$ is the required pair with $g$ equalizing $f$ and $f'$.

Thus, a functor from $\cat C^{\rm op}\lora \cat{Ab}$ is left exact if and only if it has the form as in \rref{eq3.2}. The category of left exact functor $\fun{Lex}(\cat C^\op,\cat{Ab})$ is naturally equivalent with $\cat{Ind-C}$ (cf \cite[Appendix~4.5]{artin}). We know that the category $\fun{Lex}(\cat C^\op,\cat{Ab})$ is a Grothendieck category, i.e. complete, cocomplete with a generator and filtering limits preserve exact sequences. In such a category injective envelopes exist and an injective cogenerator exists (cf. \cite[Chapter~II]{mitchell} or \cite[Chapter~V,~X]{stenstrom}).

\subsection{Extension of monoidal structures}\label{sec2.4}
Assume now that $\cat C$ is a monoidal category. Thus, we have a bifunctor $\odot:\cat C\times\cat C\lora \cat C$, which induces a bifunctor $\odot:\cat{ind-C}\times \cat{ind-C}\lora \cat{ind-C}$. Explicitly, for ind-objects $\fun X=\{X_i\}$ and $\fun Y=\{Y_j\}$, we set
\bba\label{eq3.4}\fun  X\odot\fun  Y:=\lim_{\stackrel\lora  i}\lim_{\stackrel\lora  j} X_i\odot Y_j.\eea
It is easy to see that this functor define a monoidal structure on $\cat{Ind-C}$, with the unit object being the unit object in $\cat C$. In fact, we have, for any ind-object $\fun X$:
\bbs \fun X\odot I=\lim_{\stackrel\lora  i}X_i\odot I=\lim_{\stackrel\lora  i}X_i=\fun X.\ees

Assume now that $\cat C$ is abelian.
Since directed limits preserve exact sequences, the tensor product in $\cat{Ind-C}$ is left (right) exact whenever the tensor product in $\cat C$ is. Indeed, by Lemma \ref{lem4}, any left exact sequence $0\lora \fun X\lora \fun X'\lora \fun X''$ can be ``uniformly'' represented by a filtering system
\bbs 0\lora X_i\lora X_i'\lora X''_i,\quad i\in \cat I.\ees
Thus, assuming that the tensor product on $\cat C$ is left exact, for any object $Y_j$ of $\cat C$, the sequence
\bbs 0\lora X_i\odot Y_j\lora X_i'\odot Y_j\lora X''_i\odot Y_j,\quad i\in \cat I, j\in \cat J,\ees
is exact. Since filtering limits preserve exact sequences, we have a left exact sequence
\bbs 0\lora \cat X\odot Y_j\lora\cat X'\odot Y_j\lora\cat X''\odot Y_j, j\in \cat J.\ees
Taking the limit after $j$, we obtain a left exact sequence
\bbs 0\lora\cat X\odot\cat Y\lora\cat X'\odot\cat Y\lora\cat X''\odot\cat Y.\ees

In particular, the tensor product on $\cat{Ind-C}$ is exact if the tensor product on $\cat C$ is.

\subsection{An embedding theorem for small abelian monoidal categories with exact tensor product}\label{sec2.5}
Let $\cat C$ be a small, abelian monoidal category with an exact tensor product. Set $\cat A:=\cat{Ind-C}$. Then we see in the previous section that $\cat A$ is a monoidal Grothendieck category with an exact tensor product. Let $J$ be an injective cogenerator in $\cat A$, which exists due to the the fact that $\cat A$ is a Grothendieck category.

Let $R:=\End(J)$. Consider the functor
\bbs \Hom(-,J):\cat A^{\rm op}\lora {}_R\Mod.\ees
Since $J$ is injective, the functor is exact and since $J$ is a cogenerator, the functor is faithful. Moreover, the functor is full whenever $X$ is a submodule of $J^{\oplus n}$, $n<\infty$ (cf. \cite[IV.4.1]{mitchell}). The object $J$ can be chosen so that every object of $\cat C$ satisfies this condition. Therefore, the embedding $\cat C^{\rm op}\lora {}_R\Mod$ is exact and full.

The tensor product on $\cat A$ induces a bifunctor on a subcategory of $_R\Mod$, which contains $R$ -- a progenerator of $_R\Mod$. Notice that $\Hom_{\cat A^{\rm op}}(J,J)={}_R\Hom(R,R)$. Therefore we can extend the tensor product, which is considered as a functor on the full subcategory of ${}_R\Mod_R$, consisting of one object $R$, to a colimit preserving functor on the whole category $_R\Mod$ (cf. \cite[V.5.2,p106]{mitchell}). The explicit construction is given as follows.

First, we define the tensor product on the direct sum of $R$. For any sets $S,T$, $R^S\bigcirc R^T=(R\odot R)^{S\times T}$. Then, for any module $M$, take a resolution
\bbs R^S\stackrel f\lora R^T\stackrel g\lora M\lora 0\ees
and define $R^U\bigcirc M$ to be the cokern of $R^U\bigcirc f$:
\bbs R^U\bigcirc R^S\lora R^U\bigcirc  R^T\lora R^U\bigcirc M\lora 0\ees
Analogously, we define $M\bigcirc R^U$ and then $M\bigcirc N$. Lemma V.5.2.1 of \cite{mitchell} ensures that the above construction does not depend on the choice of resolution. The associator is defined first on the direct sums of $R$ and then projected on the other objects.

From the construction of $\bigcirc$, we see that if $M$ and $N$ are finitely presented modules then
\bbs M\odot N\cong M\bigcirc N.\ees
On the other hand, we know that if an object $X$ of $\cat C$ has a resolution of the form $0\lora X\lora J^{S}\lora J^{T}$ where $S,T$ are finite sets, then the $R$-module $M=\Hom(X,J)$ is a finitely presented $R$-module.

This condition may not be satisfied for any injective cogenerator. However, it can be archived by increasing the cogenerator. We take the direct sum of all objects from $\cat C$ and then take its injective envelope. Denote the object obtained by $J_1$. Then  $J\oplus J_1$ is also an injective cogenerator, in which every object of $\cat C$ can be embedded. For any $X\in\cat C$, let $i_X$ be an embedding in $J\oplus J_1$ and let $X'$ be the cokern of $i_X$, i.e., we have an exact sequence $0\lora X\lora J\oplus J_1\lora X'\lora 0$. Now, let $J_2$ be the injective envelope of the direct sum of all $X'$ where $X$ runs in $\cat C$. Let $\bar J:=J\oplus J_1\oplus J_2$. Then for any $X\in \cat C$, we have a resolution by $\bar J$:
\bbs 
0\rTo X \rInto \bar J\rOnto J\oplus X'\oplus J_2 \rInto \bar J\oplus \bar J.\ees

Since $\bigcirc$ preserves colimit and since any module is a filtering direct limit of finite presented modules (cf. \cite[I.5]{stenstrom}, we have, for any $R$-module $\displaystyle M=\lim_{\stackrel\lora  i} M_i$, $M_i$ are finitely presented $R$-modules,
\bbs I\bigcirc M\cong I\bigcirc \lim_{\stackrel\lora  i} M_i\cong \lim_{\stackrel\lora  i} (I\bigcirc M_i)\cong \lim_{\stackrel\lora  i} M_i=M.\ees
Thus, $I$ is the unit object in $_R\Mod$ with respect to the tensor product $\bigcirc$.

Applying the result of the previous section, we have a monoidal functor 
\bbas \omega:{}_R\Mod\lora {}_R\Mod_R,\quad M\loma M\bigcirc R,\eeas
which is a right exact embedding. Compose $\omega$ with the functor $\Hom(-,\bar J)$, we get a right exact functor from $\cat A^{\rm op}$ to ${}_R\Mod_R$, whose restriction on $\cat C^{\rm op}$ is a right exact monoidal embedding. The last functor is given by
\bbs X\loma \Hom(X\odot J,J).\ees

If we start instead with $\cat C^{\rm op}$ then the above discussion give us a monoidal right exact embedding from $\cat C$ to $_R\Mod_R$. Thus we have proved
\begin{thm}\label{thm2} Let $\cat C$ be a small abelian monoidal category with the tensor product being exact. Then $\cat C$ admits a right exact monoidal embedding into the category $_R\Mod_R$ for some ring $R$. The functor is given explicitly by $X\lora \Hom(J,X\odot J)$ for a suitably chosen injective cogenerator $J$ in $\cat{Ind-C}$.\end{thm}

In particular, if $\cat C$ is an abelian rigid monoidal category then the tensor product is exact and the theorem above applies. In this case, the embedding is exact. Indeed, let $0\lora X\lora Y\lora Z\lora 0$ be an exact sequence. Applying the left dual functor, we have an exact sequence $0\lora Z^*\lora Y^*\lora X^*\lora 0$. Let $\fun F$ denote the embedding functor. Then we have a right exact sequence of $R$-modules:
\bba\label{eq3.5} \fun F(Z^*)\lora \fun F(Y^*)\lora \fun F(X^*)\lora 0\eea
Since $\fun F$ is a monoidal functor, $\fun F(X^*)\cong \fun F(X)^*=\lhom(\fun F(X),R)$. Therefore $\fun F(X)\cong \rhom(\fun F(X^*),R)$. Applying the left exact functor $\Hom(-,R)$, we obtain a left exact sequence
\bbs 0\lora \fun F(X)\lora \fun F(Y)\lora \fun F(Z).\ees

Since the category of bimodules over a ring is cocomplete, the embedding can be extended to a functor $\fun{ind-F}:\cat{ind-A}\lora {}_R\Mod_R$, which is also exact and monoidal. Explicitly, this functor has the form, for $\fun X=\{X_i\}_{i\in\fun I}$:
\bba\label{eq3.6}\fun{ind-F}(\fun X)=\lim_{\stackrel\lora i}\Hom(J,X_i\odot J).\eea

\begin{thm}\label{thm3} Let $\cat C$ be a small abelian monoidal rigid category. Then $\cat C$ admits an exact monoidal embedding in to the category of bimodules over a ring. Further, the embedding is extendable to an exact embedding of the category $\cat{ind-C}$, which commutes with colimits.\end{thm}
\proof What remains to be proved is that the functor $\fun{ind-F}$ is faithful or equivalently that $\fun{ind-F}(\fun{X})\neq 0$ when ever $\fun{X}\neq 0$. 

First, we remark that, since $\fun{\ind-F}$ is exact, $\fun{\ind-F}(\fun X)\not\cong 0$ whenever $\fun X$ possesses a subobject (or a quotient object) $\fun Y$, with $\fun{ind-F}(\fun Y)\not\cong 0$.

The following fact in $\cat{Ind-C}$ is well-known (cf. \cite[Cor.~II.3.2]{gabriel}). If and ind-object $\fun X$ is a subobject of an object $X\in\cat C$, then $\fun X$ contains a subobject $Y\in\cat C$. Indeed, let $i:\cat X\lora X$ be a monomorphism in $\cat{Ind-C}$ and $j:Y\lora\fun X$ be a non-zero morphism, $Y\in\cat C$, then $i\circ j:Y\lora X$ is non-zero and is a morphism in $\cat C$, for $\cat C$ is a full subcategory of $\cat{Ind-C}$.

A direct consequence of this fact and the preceding remark is that the image under $\fun{ind-F}$ of any non-zero ind-object, which is a subobject of an object from $\cat C$, is non-zero. 

Let now $\fun X\in\cat{Ind-C}$ be a non-zero object. There exists a non-zero morphism $j:Y\lora \lhom(\fun X,I)$, $Y\in\cat C$. Since 
\bbs \Hom(\fun X,{}^*Y)\cong\Hom(Y\odot\fun X,I)\cong\Hom(Y,\lhom(\fun X,I)),\ees
there exists a non-zero morphism $k:\fun X\lora {}^*Y$, corresponding to $j$ in the above isomorphisms. $\Im k$ is a non-zero subobject of ${}^*Y$, hence $\fun{ind-F}(\Im k)\not\cong 0$, consequently, $\fun{ind-F}(\fun X)\not\cong 0$.
\eee

\section{Semisimple abelian rigid monoidal categories}\label{sec4}

In this section we consider a simple case, where the construction in the previous section can be explicitly given. Let $\cat C$ be a semisimple rigid monoidal category. Thus, as abelian category, $\cat C$ is characterized by its simple objects $X_i,i\in \cal I$, and the ring $R_i=\End(X_i)$, where $\cal I$ is a set. The category $\cat A=\cat{Ind-C}$ is easy to characterize. Each object of $\cat A$ is a direct sum of copies of $X_i,i\in \cal I$. An injective cogenerator can be chosen to be $J=\bigoplus_{i\in \cal I}X_i$. Our embedding is then
\bba\label{eq3.7} X\loma \Hom(\bigoplus_jX_j,X_i\odot(\bigoplus_kX_k))\cong \prod_j\bigoplus_k\Hom(X_j,X_i\odot X_k),\eea
in the category of $R-R$-bimodules, where $R=\End(J)\cong\prod_iR_i$. Each $R_i$ being an endomorphism ring of a simple object is a skew-field (non-commutative field).

For any object $X\in\cat C$, we have $\End(X)\cong \Hom(I,X\ot X^*)$. If $I$ is a simple object, then the dimension of $\End(X)$ over $\bK=\End(I)$ is equal to the number of copies of $I$ in the decomposition of $X\ot X^*$. Through the above embedding, $I$ is mapped to the ring $R$ and $X\ot X^*$ is mapped to a bimodule over $R$ which is projective of finite rank when considered as left or right $R$-module. Therefore $X\ot X^*$ can contain only finitely many copies of $I$ in its decomposition into simple objects. Consequently, $\End(X)$ is finite dimension over $\bK$. 
\begin{pro} Let $\cat C$ be a semisimple rigid monoidal category with a simple unit object. Then for any object $X\in\cat C$, $\End(X)$ is finite dimensional over $\End(I)$, in particular, objects of $X$ are direct sums of simple objects.\end{pro}

From now on we shall assume that$\cat C$ is a semisimple rigid monoidal category with a simple unit object.
Let $\V ijk:=\Hom(X_j,X_i\odot X_k)$ then $\V ijk$ is an $R-R$-bimodule, whose left action is induced by the left action of $R_k$ and whose right action is induced by the left action of $R_j$. Since $R_j$ and $R_k$ are skew-fields, $\V ijk$ becomes a left and right vector space over $R_k$ and $R_j$, respectively.

Let $c_{ik}^j$ be the multiplicity of $X_j$ in $X_i\odot X_k$. Then $\dim_{R_j}\V ijk=c_{ik}^j$. On the other hand, since $\cat C$ is rigid, we have
\bba\label{eq3.8} \Hom(X_j,X_i\odot X_k)\cong\Hom(X_i{}^*\odot X_j,X_k).\eea
Therefore $\dim_{R_k} \V ijk=c_{i^*j}^k$. In particular, $\V ijk$ is finite dimensional over $R_j$ and $R_k$. Moreover, if we fix $i,j$ and let $k$ run in $\cal I$ then \rref{eq3.8} shows that  there are only finitely many $k$ for which $\V ijk$ is non-zero. Analogously, for $i,k$ fixed, there are only finitely many $j$ for which $\V ijk$ is non-zero. We also notice that $\V ijk\ot  \V lmn=0$ unless $k=n$.

Since $X_i$ is rigid, its image in ${}_R\Mod_R$ is projective of finite rank when considered as a left of a right $R$-module, thus, particularly finitely presented. Since the tensor product with a finitely presented module commute with direct product, (cf. \cite[Lemma~I.13.2]{stenstrom}), we have
\bbas \left(\prod_j\bigoplus_k\V ijk\right)\ot \left(\prod_m\bigoplus_n\V lmn\right)&\cong& 
\prod_j\bigoplus_k\left(\V ijk\ot (\prod_m\bigoplus_n\V lmn)\right)\\
&=& \prod_j\bigoplus_k\left(\V ijk\ot(\bigoplus_n\V lmn)\right)\\
&=&\prod_j\bigoplus_n\left(\bigoplus_k \V ijk\ot \V lkn\right).\eeas
Notice that in the last term of the above equations, the index $k$ runs in a finite set for each fixed $j$ and $n$. Thus, we can interpret the bimodule image of any object of $\cat C$ as an infinite matrix in which there is only finite number of non-zero element on each row or each column. The tensor product is given in terms of matrix multiplication. Dual objects correspond to transposed matrices.

The discussion above allows us to give an estimation on the dimension of the space $\End(X^{\odot n}{})$ over $\bK=\End(I)$ assuming that $I$ is simple and $R_i=\bK$ for all $i\in\cal I$. Notice that although the two actions of $\bK$ on $\Hom(X,Y)$ may be different, the dimension of $\Hom(X,Y)$ over $\bK$ with respect to these actions are equal (whenever they are finite). Thus, it is meaningful to speak of the dimension over $\bK$. In our case, $\End(X^{\odot n}{})\cong \Hom(I,X^{\odot n}{}\odot X^{\odot n}{}^*)$. Since $I$ is simple, the dimension of $\End(X^{\odot n}{})$ is equal to the number of copies of $I$ in the decomposition of $\Hom(I,X^{\odot n}\odot X^{\odot n}{}^*)$. We want to show that this dimension does not exceed $d^n$ for some positive $d$ depending only on $V$.
Embedding $\cat C$ into ${}_R\Mod_R$ as above, we see that this dimension can not exceed the number of copies of $R$ in the image of $X^{\odot n}{}\odot X^{\odot n}{}^*$. Let $V$ be the matrix representing the image of $X$, then $V^n\cdot {}^t(V^n)$ represents the image of $X^{\odot n}{}\odot X^{\odot n}{}^*$. $R$ itself is represented by the identity matrix. Therefore, the number of copies of $R$ in the bimodule $V^n\cdot {}^t(V^n)$ is equal to the minimal among the dimension over $\bK$ of the sub-bimodules lying in the diagonal of $V^n\cdot{}^t(V^n)$.

Since $X$ is rigid, $V$ represents a projective  module of finite rank over $R$, therefore, the sum of dimension of $V^j_k$ on each row or each column should be uniquely bounded by a certain number $d$. Then, the same holds for the matrix $V^n\cdot{}^t(V^n)$ with $d$ replaced by $d^{2n}$. In particular, the dimension of $(V^n\cdot{}^t(V^n))_i^i$ should not exceed $d^{2n}$. Thus, we have proved
\begin{thm}\label{thm4}Let $\cat C$ be a semisimple abelian rigid monoidal category with simple unit object, whose endomorphism ring is denoted by $\bK$. Assume that for any simple object, its endomorphism ring is isomorphic to $\bK$. Then for any object $X$, there exists a positive number $d$, such that the dimension over $\bK$ of $\End(X^n)$ does not exceed $d^n$.\end{thm}

\noindent{\it Remark.} The condition $R_i=\bK,\forall i\in\cal I$ can be replaced by the condition that the dimension of $R_i$ over $\bK$ is globally bounded by a number $c$. In this case, $d$ should be replaced by $dc^2$.

Theorem \ref{thm4} has the following important consequence.
\begin{cor}\label{cor2} Assume that $\cat C$ satisfies the condition of Theorem \ref{thm4} and that, moreover $\cat C$ is symmetric. Then, if char $\bK=0$, we can modify the symmetry on $\cat C$ so that for any object $X$ of $\cat C$, there exists an integer $n$, for which $\bigwedge_n(X)$ -- the $n$-th antisymmetric tensor power of $X$ is zero. Consequently $\cat C$ is Tannakian.\end{cor}
We give here only a sketch of the proof. A detailed proof will be given elsewhere.

Given a symmetry of $\cat C$, we define for each object $X$ its categorical dimension to be the morphism $I\stackrel\db{\textstyle\lora} X\odot X^*\stackrel\tau{\textstyle\lora} X^*\odot X\stackrel\ev{\textstyle\lora} I$, an element of $\bK=\End(I)$. This dimension if an additive and tensor-multiplicative function on $X$. Since the category is semisimple, the dimension of a simple object is non-zero.

On the other hand, the symmetry induces a representation of $k[\Ss_n]$ in $\End(X^{\odot n})$ for any object $X$, $\Ss_n$ is the symmetric group. Theorem \ref{thm4} ensures that starting form some $n$, the representation is not faithful. That means some subobject of $X^{\odot n}$ should be zero; its dimension is therefore also zero. This implies that the categorical dimension of $X$ should be an integer. Thanks the semisimplicity we can modify the symmetry on $\cat C$ so that the dimension of a simple object is a positive integer, hence so is the dimension of any object. Then we are done by \cite[Theorem~7.1]{deligne90}.\eee

\begin{center}\bf Acknowledgment\end{center}
\nopagebreak
\mythanks


\begin{thebibliography}{10}

\bibitem{sga4}
M.~Artin, A.~Grothendieck, and J.~L. Verdier.
\newblock {\em Th\'eorie des topos et cohomologie \'etale des schémas.
  S\'eminaire de G\'eom\'etrie Alg\'ebrique du Bois-Marie 1963--1964 (SGA 4)},
  volume 305 of {\em Lecture Notes in Mathematics}.
\newblock Springer-Verlag, 1973.

\bibitem{artin}
M.~Artin and B~Mazur.
\newblock {\em Etale homotopy}, volume No. 100 of {\em Lecture Notes in
  Mathematics}.
\newblock Springer-Verlag, Berlin-New York, 1969.

\bibitem{deligne90}
P~Deligne.
\newblock Cat\'egories tannakiennes.
\newblock In Cartier P. and et.al., editors, {\em The Grothendieck
  Festschrift}, volume~II of {\em Progr. Math., 87}, pages 111--195.
  Birkh\"auser Boston, Boston, MA, 1990.

\bibitem{deligne82}
P.~{D}eligne and J.~{M}ilne.
\newblock {Tannakian Categories}.
\newblock In {\em Lecture Notes in Mathematics}, volume 900, pages 101--228.
  Springer-Verlag, Berlin-Heidelberg-New York, 1982.

\bibitem{dr89}
Sergio Doplicher and John~E. Roberts.
\newblock A new duality theory for compact groups.
\newblock {\em Invent. Math.}, 98(1):157--218, 1989.

\bibitem{drinfeld86}
V.G. {D}rinfel'd.
\newblock {Quantum Groups}.
\newblock {\em Proceedings of ICM, Berkeley}, 1987.

\bibitem{gabriel}
P.~Gabriel.
\newblock Des cat\'egories ab\'eliennes.
\newblock {\em Bull. Soc. Math. France}, 90:323--448, 1962.

\bibitem{kassel}
Ch. Kassel.
\newblock {\em Quantum Groups}, volume 155 of {\em Graduate Texts in
  Mathematics}.
\newblock Springer-Verlag, 1995.
\newblock 531p.

\bibitem{maclane}
S.~Mac Lane.
\newblock {\em Categories. {F}or the Working Mathematician}.
\newblock Springer Verlag, 1971.

\bibitem{mitchell}
B.~Mitchell.
\newblock {\em Theory of categories}, volume XVII of {\em Pure and Applied
  Mathematics}.
\newblock Academic Press, New York-London, 1965.

\bibitem{par3}
B.~Pareigis.
\newblock Lectures on quantum groups.
\newblock Available at www.mathematik.uni-muenchen.de/\verb-~-pareigis, 1994.

\bibitem{schauen1}
P.~Schauenburg.
\newblock The monoidal center construction and bimodules.
\newblock {\em To appear in Journal of Pure and Applied Algebra}, 2000.

\bibitem{stenstrom}
B.~Stenstr\"om.
\newblock {\em Rings of quotients}, volume 217 of {\em Die Grundlehren der
  Mathematischen Wissenschaften}.
\newblock Springer-Verlag, New York-Heidelberg, 1975.

\bibitem{watts}
C.~E. Watts.
\newblock Intrinsic characterizations of some additive functors.
\newblock {\em Proc. Amer. Math. Soc.}, 11:5--8, 1960.

\end{thebibliography}
\end{document}